\documentclass[onecolumn,11pt,draftcls]{IEEEtran}
\usepackage{amsmath,amssymb,amsthm,graphicx,graphics,subfigure,epsfig,enumerate}
\usepackage{parskip}		
\usepackage{bm}
\usepackage{color}

\usepackage{blkarray}
\usepackage{multirow}
\usepackage{mathtools}
\usepackage{url}  
\usepackage{ifthen}  
\usepackage[utf8]{inputenc}

\newcommand\remark[1]{\textbf{Remark #1:}}

\newcommand{\gbh}{{\hat{\bf g}}}
\newcommand{\Gb}{{\bf G}}
\newcommand{\Gbh}{{\hat{\bf G}}}
\newcommand{\Jb}{{\bf J}}
\newcommand{\Jbh}{{\hat{\bf J}}}

\newcommand{\thetab}{{\boldsymbol{\theta}}}
\newcommand{\eb}{{\bf e}}
\newcommand{\xb}{{\bf x}}

\newcommand{\Vb}{{\bf V}}

\newcommand{\Phib}{{\bf \Phi}}

\newcommand{\Wb}{{\bf W}}

\begin{document}

\title{Analysis of Fisher Information and the Cram\'{e}r-Rao Bound for Nonlinear Parameter Estimation after Compressed Sensing}
\author{Pooria Pakrooh,~\IEEEmembership{Student Member,~IEEE,} Ali~Pezeshki,~\IEEEmembership{Member,~IEEE,} Louis L. Scharf,~\IEEEmembership{Life Fellow,~IEEE,} Douglas Cochran,~\IEEEmembership{Senior Member,~IEEE,} and Stephen D. Howard\vfill
\thanks{A preliminary version of parts of this paper was presented at the 2013 IEEE International Conference on Acoustics, Speech, and Signal Processing (ICASSP'13), Vancouver, BC, Canada, May 26-31, 2013.}
\thanks{This work is supported in part by NSF under grants CCF-1018472 and CCF-1422658.}
\thanks{P. Pakrooh is with the Department of Electrical and
Computer Engineering, Colorado State University, Fort Collins, CO 80523, USA (e-mail: Pooria.Pakrooh@colostate.edu).}
\thanks{A. Pezeshki is with the Department of Electrical and
Computer Engineering and the Department of Mathematics, Colorado State University, Fort Collins, CO 80523, USA (e-mail: Ali.Pezeshki@colostate.edu).}
\thanks{L. L. Scharf is with the Department of Mathematics and the Department of Statistics, Colorado State University, Fort Collins, CO 80523, USA (e-mail: Louis.Scharf@colostate.edu).}
\thanks{D. Cochran is with the School of Electrical, Computer, and Energy Engineering, Arizona State University, Tempe, AZ, 85287, USA (e-mail: cochran@asu.edu).}
\thanks{S. D. Howard is with the Defence Science and Technology Organisation, Edinburgh, SA 5111, Australia (e-mail: sdhoward@unimelb.edu.au).}}

\maketitle

\begin{abstract}

In this paper, we analyze the impact of compressed sensing with complex random matrices on Fisher information and the Cram\'{e}r-Rao Bound (CRB) for estimating unknown parameters in the mean value function of a complex multivariate normal distribution. We consider the class of random compression matrices whose distribution is right-orthogonally invariant. The compression matrix whose elements are i.i.d. standard normal random variables is one such matrix. We show that for all such compression matrices, the Fisher information matrix has a complex matrix  beta distribution. We also derive the distribution of CRB. These distributions can be used to quantify the loss in CRB as a function of the Fisher information of the non-compressed data. In our numerical examples, we consider a direction of arrival estimation problem and discuss the use of these distributions as guidelines for choosing compression ratios based on the resulting loss in CRB. 

\end{abstract}

\begin{keywords}
Cram\'{e}r-Rao bound, compressed sensing, Fisher information, matrix beta distribution, parameter estimation \vspace{-.2cm}
\end{keywords}

\section{Introduction} \label{sec:intro}
Inversion of a measurement for its underlying modes is an important topic which has applications in communications, radar/sonar signal processing and optical imaging. Classical methods for inversion are based on maximum likelihood, variations on linear prediction, subspace filtering, etc. Compressed sensing \cite{CS1}--\nocite{CS2}\cite{CS3} is a relatively new theory which exploits sparse representations and sparse recovery for inversion.

In our previous work \cite{chi-sp11}--\nocite{chi-icassp10,scharf-asilomar11}\cite{scharf-DASP11}, the sensitivity of sparse inversion algorithms to basis mismatch and frame mismatch were studied. Our results show that mismatch between the actual basis in which a signal has a sparse representation and the basis (or frame) which is used for sparsity in a sparse reconstruction algorithm, e.g., basis pursuit, has performance consequences on the estimated parameter vector.

This paper addresses another fundamental question: \emph{How much information is  retained (or lost) in compressed noisy measurements for nonlinear parameter estimation?} To answer this question, we analyze the effect of compressed sensing on the Fisher information matrix and the Cram\'{e}r-Rao bound (CRB). We derive the distribution of the Fisher information matrix and the CRB for the class of random matrices whose distributions are invariant under right-orthogonal transformations. These include i.i.d draws of spherically invariant matrix rows, including, for example, i.i.d. draws of standard normal matrix elements.

Our prior work on compressed sensing and the Fisher information matrix \cite{scharf-asilomar11,scharf-DASP11} contain numerical results that characterize the increase in CRB after random compression for the case where the parameters nonlinearly modulate the mean in a multivariate normal measurement model. Also, in \cite{Pakrooh13a} we derived analytical lower and upper bounds on the CRB for the same problem, and used our bounds to quantify the potential loss in CRB.

Other studies on the effect of compressed sensing on the CRB and the Fisher information matrix include \cite{Babadi09}--\nocite{Babaizadeh12}\cite{Madhow12}. Babadi \textit{et al.} \cite{Babadi09} proposed a ``Joint Typicality Estimator'' to show the existence of an estimator that asymptotically achieves the CRB of sparse parameter estimation for random Gaussian compression matrices. Niazadeh \textit{el al.} \cite{Babaizadeh12} generalize the results of \cite{Babadi09} to a class of random compression matrices which satisfy the concentration of measures inequality. Nielsen \textit{et al.} \cite{Nielsen2012} derive the mean of the Fisher information matrix for the same class of random compression matrices that we are considering. Ramasamy \textit{et al.} \cite{Madhow12}  derive bounds on the Fisher information matrix, but not for the model we are considering. We will clarify the distinction between our work and \cite{Madhow12} after establishing our notation in Section \ref{sec:prob}.

\section{Problem statement} \label{sec:prob}
Let $\mathbf{y}\in \mathbb{C}^n$ be a complex random vector whose probability density function $f(\mathbf{y};\boldsymbol{\theta})$ is parameterized by an unknown but deterministic parameter vector $\boldsymbol{\theta}\in \mathbb{R}^p$. The derivative of the log-likelihood function with respect to $\boldsymbol{\theta}=[\theta_1,\theta_2, \cdots,\theta_p]$ is called the Fisher score, and the covariance matrix of the Fisher score is the Fisher information matrix which we denote by $\mathbf{J(\boldsymbol{\theta})}$:
\begin{equation}\label{eq:Fisher}
\mathbf{J(\boldsymbol{\theta})}=E[(\frac{\partial \log f(\mathbf{y};\boldsymbol{\theta})}{\partial \boldsymbol{\theta}})(\frac{\partial \log f(\mathbf{y};\boldsymbol{\theta})}{\partial \boldsymbol{\theta}})^H]\in\mathbb{C}^{p\times p}.
\end{equation}
The inverse $\mathbf{J}^{-1}(\boldsymbol{\theta})$ of the Fisher information matrix lower bounds the error covariance matrix for any unbiased estimator $\hat{\boldsymbol{\theta}}(\mathbf{y})$ of $\boldsymbol{\theta}$, that is
\begin{equation}
E[(\hat{\boldsymbol{\theta}}(\mathbf{y})-\boldsymbol{\theta})(\hat{\boldsymbol{\theta}}(\mathbf{y})-\boldsymbol{\theta})^H]\succeq \mathbf{J}^{-1}(\boldsymbol{\theta}),
\end{equation}
where $\mathbf{A}\succeq \mathbf{B}$ for matrices $\mathbf{A},\mathbf{B}\in\mathbb{C}^{n\times n}$ means $\mathbf{a}^H\mathbf{A}\mathbf{a}\geq \mathbf{a}^H\mathbf{B}\mathbf{a}$ for all $\mathbf{a}\in \mathbb{C}^n$. The $i^{th}$ diagonal element of $\mathbf{J}^{-1}(\boldsymbol{\theta})$ is the CRB for estimating $\theta_i$ and it gives a lower bound on the MSE of any unbiased estimator of $\theta_i$ from $\mathbf{y}$ (see, e.g., \cite{4-Eldar}).

For $\mathbf{y}\in \mathbb{C}^n$ a proper complex normal random vector distributed as $\mathcal{CN}(\mathbf{x}(\boldsymbol{\theta}),\mathbf{C})$ with unknown mean vector $\mathbf{x}(\boldsymbol{\theta})$ parameterized by $\boldsymbol{\theta}$, and known covariance $\mathbf{C}=\sigma^2 \mathbf{I}$, the Fisher information matrix is the Grammian
\begin{equation}\label{eq:FisherGaussian}
\mathbf{J(\boldsymbol{\theta})}= \mathbf{G}^H\mathbf{C}^{-1}\mathbf{G}=\frac{1}{\sigma^2}\mathbf{G}^H\mathbf{G}.
\end{equation}
The $i^{th}$ column $\mathbf{g}_i$ of $\mathbf{G}=[\mathbf{g}_1, \mathbf{g}_2, \cdots,\mathbf{g}_p]$ is the partial derivative $\mathbf{g}_i=\frac{\partial}{\partial{\theta}_i}\mathbf{x}(\boldsymbol{\theta})$, which characterizes the sensitivity of the mean vector $\mathbf{x}(\boldsymbol{\theta})$ to variation of the $i^{th}$ parameter $\mathbf{\theta}_i$. 

The CRB for estimating $\theta_i$ is given by
\begin{equation} {(\mathbf{J}^{-1}(\boldsymbol{\theta}))}_{ii}=\sigma^2(\mathbf{g}_i^H(\mathbf{I}-\mathbf{P_{\mathbf{G}_i}})\mathbf{g}_i)^{-1},
\end{equation}
where $\mathbf{G}_i$ consists of all columns of $\mathbf{G}$ except $\mathbf{g}_i$, and $\mathbf{P}_{\mathbf{G}_i}$ is the orthogonal projection onto the column space of $\mathbf{G}_i$ \cite{Scharf_geo}. This CRB can also be written as
\begin{equation}
{(\mathbf{J}^{-1}(\boldsymbol{\theta}))}_{ii}=\frac{\sigma^2}{\|\mathbf{g}_i\|_2^2\sin^2(\psi_i)},
\end{equation}
where $\psi_i$ is the principal angle between subspaces $\langle\mathbf{g}_i\rangle$ and $\langle\mathbf{G}_i\rangle$. These representations illuminate the geometry of the CRB, which is discussed in detail in \cite{Scharf_geo}.

If $\mathbf{y}$ is compressed by a compression matrix $\mathbf{\Phi}\in\mathbb{C}^{m\times n}$ to produce $\hat{\mathbf{y}}=\mathbf{\Phi}\mathbf{y}$, then the probability density function of the compressed data $\hat{\mathbf{y}}$ is $\mathcal{CN}[\mathbf{\Phi}\mathbf{x}(\mathbf{\theta}),\sigma^2\mathbf{\Phi}\mathbf{\Phi}^H]$. The Fisher information matrix $\hat{\mathbf{J}}(\boldsymbol{\theta})$ is given by
\begin{equation}\label{eq:FisherAfter}
\hat{\mathbf{J}}(\boldsymbol{\theta})=\frac{1}{\sigma^2}\hat{\mathbf{G}}^H\hat{\mathbf{G}},
\end{equation}
where $\hat{\mathbf{G}}=\mathbf{P}_{\mathbf{\Phi}^H}\mathbf{G}$. The CRB for estimating the $i^{th}$ parameter is
\begin{equation} {(\Jbh^{-1}(\thetab))}_{ii}=\sigma^2(\gbh_i^H(\mathbf{I}-\mathbf{P_{\Gbh_i}})\gbh_i)^{-1},
\end{equation}
where $\hat{\mathbf{G}}_i=\mathbf{P}_{\mathbf{\Phi}^H}\mathbf{G}_i$, and $\mathbf{P}_{\mathbf{\Phi}^H}=\mathbf{\Phi}^H(\mathbf{\Phi}\mathbf{\Phi}^H)^{-1}\mathbf{\Phi}$ is the orthogonal projection onto the row span of $\mathbf{\Phi}$.

Our aim is to study the effect of random compression on the Fisher information matrix and the CRB. In section \ref{sec:Distn} we investigate this problem by deriving the distributions of the Fisher information matrix and the CRB for the case in which the elements of the compression matrix $\Phi$ are distributed as i.i.d. standard normal random variables. Then we demonstrate that the same analysis holds for a wider range of random compression matrices. 

\remark{1} In parallel to our work, Ramasamy et al. \cite{Madhow12} have also looked at the impact of compression on Fisher information. However, they have considered a different parameter model. Specifically, their compressed data has density $\mathcal{CN}[\mathbf{\Phi}\mathbf{x}(\boldsymbol{\theta}),\sigma^2\mathbf{I}]$, in contrast to ours which is distributed as $\mathcal{CN}[\mathbf{\Phi}\mathbf{x}(\boldsymbol{\theta}),\sigma^2\mathbf{\Phi}\mathbf{\Phi}^H]$. Our model is a signal-plus-noise model, wherein the noisy signal $\mathbf{x}({\boldsymbol{\theta}})+\mathbf{n}$, $\mathbf{n}\sim\mathcal{CN}(\mathbf{0},\sigma^2\mathbf{I})$, is compressed to produce $\mathbf{\Phi}\mathbf{x}(\boldsymbol{\theta})+\mathbf{\Phi}\mathbf{n}$. In contrast, their model corresponds to compressing a noiseless signal $\mathbf{x}(\boldsymbol{\theta})$ to produce $\mathbf{\Phi}\mathbf{x}(\theta)+\mathbf{w}$, where $\mathbf{w}\sim\mathcal{CN}(\mathbf{0},\sigma^2\mathbf{I})$ represents post-compression noise. Note that the Fisher information, CRB and corresponding bounds of these two models are different, as in our model noise enters at the input of the compressor, whereas in \cite{Madhow12} noise enters at the output of the compressor. This is an important distinction.

\section{Distribution of the Fisher information Matrix after compression} \label{sec:Distn}

Let $\mathbf{W}$ be the normalized Fisher information matrix after compression, defined as
\begin{equation}\label{eq:FM1}
\mathbf{W}=\mathbf{J}^{-1/2}\hat{\mathbf{J}}\mathbf{J}^{-H/2}\in \mathbb{C}^{p\times p},
\end{equation}
where $\Jb$ and $\Jbh$ are the Fisher information matrix before compression and the Fisher information matrix after compression, defined in (\ref{eq:FisherGaussian}) and (\ref{eq:FisherAfter}), respectively. Our aim is to derive the distribution of the random matrix $\mathbf{W}$ for the case that the elements of the compression matrix $\mathbf{\Phi}_{ij}$ are i.i.d. random variables distributed as $\mathcal{CN}(0,1)$. We assume $n-p\ge m$, which is typical in almost all compression scenarios of interest. 

Using (\ref{eq:FisherGaussian}) and (\ref{eq:FisherAfter}), $\mathbf{W}$ may be written as
\begin{align}\label{eq:FM2}
\mathbf{W}&=\mathbf{H}^H\mathbf{P}_{\mathbf{\Phi}^H}\mathbf{H}, 
\end{align}
where $\mathbf{H}=\Gb(\Gb^H\Gb)^{-H/2}$ is a left orthogonal matrix, i.e., $\mathbf{H}^H\mathbf{H}=\mathbf{I}_{p}$.
Define $\tilde{\mathbf{H}}\in\mathbb{C}^{n\times(n-p)}$ such that ${\bm\Lambda}=[\mathbf{H}|\tilde{\mathbf{H}}]$ is an orthonormal basis for $\mathbb{C}^n$, i.e., ${\bm\Lambda}{\bm\Lambda}^H={\bm\Lambda}^H{\bm\Lambda}=\mathbf{I}_n$. Then we have:
\begin{align}\label{eq:FM3}
\mathbf{W}&=\left[\begin{array}{cc}
\mathbf{I}_p\hspace{1mm}\vline&\hspace{-1mm}\mathbf{0} 
\end{array}\right]\mathbf{\Lambda}^H\mathbf{P}_{\mathbf{\Phi}^H}\mathbf{\Lambda}\left[
\begin{array}{c}
\mathbf{I}_p\\ \hline
\mathbf{0} 
\end{array}\right],
\end{align}
where 
\begin{align}\label{eq:FM3.5}
\mathbf{\Lambda}^H\mathbf{P}_{\mathbf{\Phi}^H}\mathbf{\Lambda}&=\mathbf{\Lambda}^H\mathbf{\Phi}^H(\mathbf{\Phi}\mathbf{\Phi}^H)^{-1}\mathbf{\Phi}\mathbf{\Lambda}\nonumber\\
&=\mathbf{\Lambda}^H\mathbf{\Phi}^H(\mathbf{\Phi}\mathbf{\Lambda}\mathbf{\Lambda}^H\mathbf{\Phi}^H)^{-1}\mathbf{\Phi}\mathbf{\Lambda}\nonumber\\
&=\mathbf{P}_{\mathbf{\Lambda}^H\mathbf{\Phi}^H}.
\end{align}
Because the distribution of $\mathbf{\Phi}$ is right orthogonally invariant, the distribution of $\mathbf{\Lambda}^H\mathbf{P}_{\mathbf{\Phi}^H}\mathbf{\Lambda}$ is the same as the distribution of $\mathbf{P}_{\mathbf{\Phi}^H}$. Therefore, the distribution of $\mathbf{W}$ is the same as the distribution of $\mathbf{V}=\Phib_1^H(\Phib\Phib^H)^{-1}\Phib_1$, where $\Phib=[\Phib_1|\Phib_2] $, $\Phib_1\in\mathbb{C}^{m\times p}$ and $\Phib_2\in\mathbb{C}^{m\times (n-p)}$. Now, write $\mathbf{V}$ as $\mathbf{V}=\mathbf{Y}\mathbf{Y}^H$, with $\mathbf{Y}=\Phib_1^H\mathbf{Z}^{-1/2}$ and $\mathbf{Z}=\mathbf{\Phi}\mathbf{\Phi}^H$. Since $\mathbf{Z}=\Phib_1\Phib_1^H+\Phib_2\Phib_2^H$, and  $\Phib_2\Phib_2^H$ is distributed as complex Wishart $\mathcal{W}_c(\mathbf{I}_m,m,n-p)$ for $n-p\geq m$, given $\Phib_1$ the pdf of $\mathbf{Z}$ is
\begin{equation}\label{eq:FM4}
f(\mathbf{Z}|\Phib_1)=c_1e^{-tr(\mathbf{Z}-\Phib_1\Phib_1^H)}|\mathbf{Z}-\Phib_1\Phib_1^H|^{n-m-p}.
\end{equation}
The pdf of $\Phib_1$ can be written as $c_2e^{-tr(\mathbf{\Phib_1\Phib_1^H})}$. Therefore, the joint pdf of $\mathbf{Z}$ and $\Phib_1$ is
\begin{equation}\label{eq:FM5}
f(\mathbf{Z},\Phib_1)=c_3e^{-tr(\mathbf{Z})}|\mathbf{Z}-\Phib_1\Phib_1^H|^{n-m-p}
\end{equation}
where $c_1$, $c_2$, and $c_3=c_1c_2$ are normalization factors. Since $\mathbf{Y}=\Phib_1^H\mathbf{Z}^{-1/2}$, from (\ref{eq:FM5}) the joint pdf of $\mathbf{Y}$ and $\mathbf{Z}$ is
\begin{align}\label{eq:FM6}
f(\mathbf{Y},\mathbf{Z})&=c_3e^{-tr(\mathbf{Z})}|\mathbf{Z}-\mathbf{Z}^{1/2}\mathbf{Y}^H\mathbf{Y}\mathbf{Z}^{H/2}|^{n-m-p}|\mathbf{Z}|^p\nonumber\\
&=c_3e^{-tr(\mathbf{Z})}|\mathbf{I}_m-\mathbf{Y}^H\mathbf{Y}|^{n-m-p}|\mathbf{Z}|^{n-m}.
\end{align}
This shows that $\mathbf{Y}$ and $\mathbf{Z}$ are independent and the pdf of $\mathbf{Y}$ is
\begin{align}\label{eq:FM7}
f(\mathbf{Y})&=c_4|\mathbf{I}_m-\mathbf{Y}^H\mathbf{Y}|^{n-m-p}\nonumber\\
&=c_4|\mathbf{I}_p-\mathbf{Y}\mathbf{Y}^H|^{n-m-p},
\end{align}
where $c_4$ is a normalization factor. 

Let $g(\mathbf{Y}\mathbf{Y}^H)=f(\mathbf{Y})$. To derive the distribution of $\mathbf{V}=\mathbf{Y}\mathbf{Y}^H$ we use the following theorem.

\textbf{Theorem 1:} \cite{Srivastava} If the density of $\mathbf{Y}\in\mathbb{C}^{p\times m}$ is $g(\mathbf{Y}\mathbf{Y}^H)$, then the density of $\mathbf{V}=\mathbf{Y}\mathbf{Y}^H$ is
 \begin{equation}\label{eq:FM8}
\frac{|\mathbf{V}|^{m-p}g(\mathbf{V})\pi^{mp}}{\tilde{\Gamma}_p(m)},
\end{equation}
where $\tilde{\Gamma}_m(.)$ is the complex multivariate Gamma function.

Using Theorem 1 and \eqref{eq:FM7}, the pdf of $\mathbf{V}$ is 
\begin{equation}\label{eq:FM9}
c_5|\mathbf{V}|^{m-p}|\mathbf{I}_p-\mathbf{V}|^{n-m-p} \hspace{4mm}\textrm{for} \hspace{4mm} \mathbf{0}\leq\mathbf{V}\leq\mathbf{I}_p,
\end{equation} 
which is the Type I complex multivariate beta distribution $\mathbb{C}B_p^I(m,n-m)$ for $c_5=\frac{\tilde{\Gamma}_p(n)}{\tilde{\Gamma}_p(m)\tilde{\Gamma}_p(n-m)}$. Recall that the distribution of the normalized Fisher information matrix after compression $\Wb=\Jb^{-1/2}{\Jbh}\Jb^{-H/2}$ is identical to that of $\Vb$. Therefore, $\Wb$ is also distributed as $\mathbb{C}B_p^I(m,n-m)$.

\remark{2} It is important to note that the distribution of  $\mathbf{W}=\Jb^{-1/2}{\Jbh}\Jb^{-H/2}$ is \emph{invariant} to $\Jb$, and it  depends on only on the parameters $(n-m)-p$ and $m-p$. In this sense, this result for the distribution of $\Jb^{-1/2}{\Jbh}\Jb^{-H/2}$ is universal, and reminiscent of the classical result of Reed, Mallat, and Brennan \cite{RMB} for normalized SNR in adaptive filtering.  

\textbf{Lemma 1:} \cite{James1964} Assume $\mathbf{A}\in\mathbb{C}^{p\times p}$ is a positive definite Hermitian random matrix with a 
pdf $h(\mathbf{A})$. Then, the joint pdf of eigenvalues $\mathbf{\Lambda}=diag(\lambda_1,\lambda_2,\dots,\lambda_p)$
of $\mathbf{A}$ is
\begin{equation}\label{lemma1}
\frac{\pi^{p(p-1)}}{\tilde{\Gamma}_p(p)}\prod_{i<j}^{p}(\lambda_i-\lambda_j)^2\int_{\mathcal{O}(p)}h(\mathbf{U}\mathbf{\Lambda}\mathbf{U}^H)d\mathbf{U},
\end{equation}
where $d\mathbf{U}$ is the invariant Haar measure on the Unitary group $\mathcal{O}(p)$.

Using Lemma 1, we can derive the joint distribution of the eigenvalues of $\Wb$. Replacing the pdf of $\Wb\sim \mathbb{C}B_p^I(m,n-m)$ in (\ref{lemma1}), the joint pdf of the eigenvalues $\lambda_1,\lambda_2,\dots,\lambda_p$ of $\mathbf{W}=\Jb^{-1/2}{\Jbh}\Jb^{-H/2}$ is given by
\begin{align}\label{eq:FM_eig}
\frac{\pi^{p(p-1)}\tilde{\Gamma}_p(n)}{\tilde{\Gamma}_p(p)\tilde{\Gamma}_p(m)\tilde{\Gamma}_p(n-m)}\prod_{i<j}^{p}(\lambda_i-\lambda_j)^2\prod_{i=1}^p\lambda_i^{m-p}(1-\lambda_i)^{n-m-p}.
\end{align}
Now, from (\ref{eq:FM9}) and using the transformation $\Jbh=\Jb^{1/2}\Wb\Jb^{H/2}$, the Fisher information matrix after compression $\Jbh$ is distributed as
\begin{equation}\label{eq:FM10}
c_5|\Jb|^{p-n}|\Jbh|^{m-p} |\Jb-\Jbh|^{n-m-p}\hspace{4mm}\textrm{for}\hspace{4mm} \mathbf{0}\leq\Jbh\leq\Jb
\end{equation} 
and the inverse of the Fisher information matrix after compression $\hat{\mathbf{K}}=\Jbh^{-1}$ is distributed as
\begin{equation}\label{eq:FM11}
c_5|\Jb|^{p-n}|\hat{\mathbf{K}}|^{-n} |\Jb\hat{\mathbf{K}}-\mathbf{I}_p|^{n-m-p}\hspace{4mm}\textrm{for}\hspace{4mm} \hat{\mathbf{K}}\geq\Jb^{-1}.
\end{equation}

\remark{3} For the class of random compression matrices that have density functions of the form $g(\Phib\Phib^H)$, that is, the distribution of $\Phib$ is right orthogonally invariant, $\Phib^H(\Phib\Phib^H)^{-1/2}$ is uniformly distributed on the Stiefel manifold $\mathcal{V}_{m}(\mathbb{C}^n)$ \cite{Chikuse2003}. Therefore, the distribution of the normalized Fisher information matrix for this class of compression matrices is the same as the one given in (\ref{eq:FM9}).

\remark{4} Using the properties of a complex multivariate beta distribution \cite{Gupta2005}, we have:
\begin{equation}
E[\Jbh]=\frac{m}{n}\Jb,
\end{equation}
and
\begin{equation}
E[\Jbh^{-1}]=\frac{n-p}{m-p}\Jb^{-1}.
\end{equation}
This shows that, on average, compression results in a factor $\frac{n-m}{n}$ loss in the Fisher information and a factor $\frac{n-p}{m-p}$ increase in the CRB $\Jb^{-1}$.

The distribution of the CRB after compression can be derived using the following Lemma.

\textbf{Lemma 2:}  \cite{Gupta2005} Assume $\mathbf{X}\sim \mathbb{C}B_p^I(a_1,a_2)$. Let $\mathbf{z}$ be a complex vector independent of $\mathbf{X}$. Then, $x=\frac{\mathbf{z}^H\mathbf{z}}{\mathbf{z}^H\mathbf{X}^{-1}\mathbf{z}}$ is distributed as $B_p^I(a_1-p+1,a_2)$, which is a Type I univariate beta distribution with the pdf
\begin{equation}
\frac{\Gamma_p(a_1+a_2)}{\Gamma_p(a_1)\Gamma_p(a_2)}x^{a_1-1}(1-x)^{a_2-1}\hspace{4mm}\textrm{for}\hspace{4mm} 0<x<1.
\end{equation}

Now consider the CRB on an unbiased estimator of parameter $\theta_i$, after compression, normalized by the CRB before compression: 
\begin{equation}\label{eq:normJ}
\frac{(\Jbh^{-1})_{ii}}{(\Jb^{-1})_{ii}}=\frac{\mathbf{e}_i^H\Jbh^{-1}\mathbf{e}_i}{\mathbf{e}_i^H\Jb^{-1}\mathbf{e}_i}
=\frac{\mathbf{z}^H\mathbf{W}^{-1}\mathbf{z}}{\mathbf{z}^H\mathbf{z}},
\end{equation} 
where $\mathbf{z}=\Jb^{-1/2}\mathbf{e}_i$, and $\mathbf{e}_i\in\mathbb{C}^{p}$ is a standard unit vector with $1$ as its $i$th element and zeros as its other elements. By Lemma 1, the above ratio is distributed as the inverse of a univariate beta random variable $B^I(m-p+1,n-m)$.  

\remark{5}  From the distribution of the CRB after compression, we have:
\begin{align}
E[(\Jbh^{-1})_{ii}]&=(\frac{n-p}{m-p})(\Jb^{-1})_{ii}, \\ \mathrm{var}[(\Jbh^{-1})_{ii}]&=\frac{(n-m)(n-p)}{(m-p-1)(m-p)^2}((\Jb^{-1})_{ii})^2.
\end{align}
\remark{6} We can also look at the effect of compressed sensing on the Kullback Leibler (KL) divergence of two normal probability distributions, for the class of random compressors already discussed in Remark 3. The KL divergence between $\mathcal{CN}(\xb(\thetab),\mathbf{C})$ and $\mathcal{CN}(\xb(\thetab '),\mathbf{C})$ is:
\begin{equation}
D(\thetab,\thetab ')=(\xb(\thetab)-\xb(\thetab '))^H\mathbf{C}^{-1}(\xb(\thetab)-\xb(\thetab ')).
\end{equation}
After compression with $\Phib$ we have 
\begin{align}
\hat{D}(\thetab,\thetab ')=(\xb(\thetab)-\xb(\thetab '))^H\Phib^H(\Phib\mathbf{C}\Phib^H)^{-1}\Phib(\xb(\thetab)-\xb(\thetab ')).
\end{align}
For the case $\mathbf{C}=\sigma^2\mathbf{I}$, the normalized KL divergence is 
\begin{equation}
\frac{\hat{D}(\thetab,\thetab ')}{D(\thetab,\thetab ')}=\frac{(\xb(\thetab)-\xb(\thetab '))^H{\mathbf{P}_{\Phib^H}} (\xb(\thetab)-\xb(\thetab '))}
{(\xb(\thetab)-\xb(\thetab '))^H(\xb(\thetab)-\xb(\thetab '))}.
\end{equation}
Therefore, the normalized KL divergence, for random compression matrices $\mathbf{\Phi}$ whose distributions are invariant to right-orthogonal transformations, is distributed as $B^I(m,n-m)$.

\section{Numerical results}

\begin{figure}[ht]\centering
  \includegraphics[width=288pt]{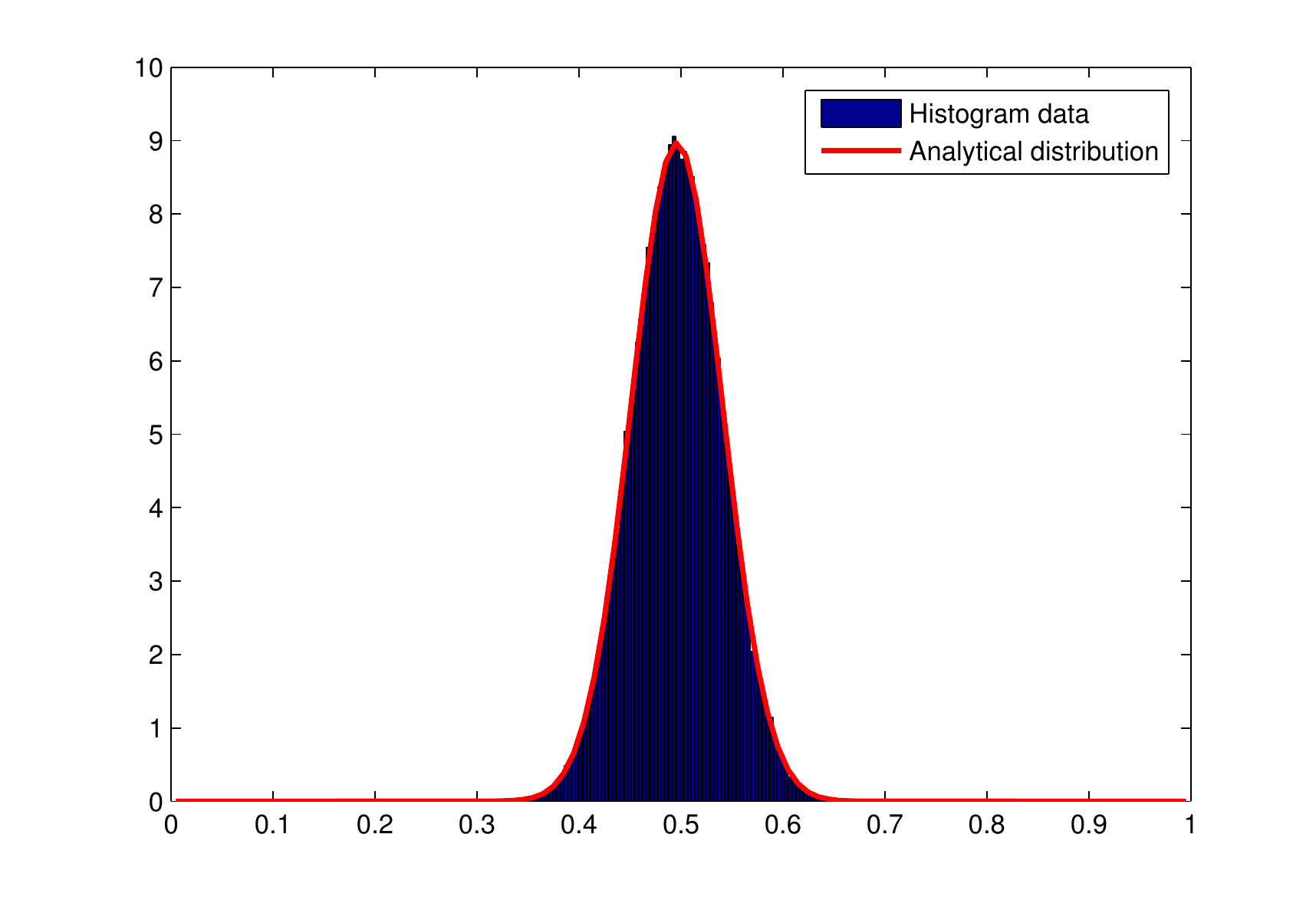}
  \caption{Histogram data and analytical distributions for $\frac{(\Jb^{-1})_{11}}{(\Jbh^{-1})_{11}}$ using $10^5$ realizations of i.i.d. Gaussian compression matrices with $n=128$ and $m=64$.}\label{fig:Hist}
\end{figure}

As a special example, we consider the effect of compression on DOA estimation using a uniform line array with $n$ elements. In our simulations, we consider two sources whose electrical angles $\mathbf{\theta}_1$ and $\mathbf{\theta}_2$ are unknown. The mean vector $\mathbf{x}(\boldsymbol{\theta})$ is $\mathbf{x}(\boldsymbol{\theta})=\mathbf{x}(\boldsymbol{\theta}_1)+\mathbf{x}(\boldsymbol{\theta}_2)$, where
\begin{equation}
\mathbf{x}(\theta_i)=A_ie^{j\phi_i}[1\hspace{2mm}e^{j\theta_i}\hspace{2mm}e^{j2\theta_i}\hspace{2mm}\cdot\cdot\cdot\hspace{2mm}e^{j(n-1)\theta_i}]^T.
\end{equation}
Here $A_i$ and $\phi_i$  are the amplitude and phase of the $i^{th}$ source, which we assume known. We set $\phi_1=\phi_2=0$ and $A_1=A_2=1$. We wish to estimate $\theta_1$, whose true value in this example is zero, in the presence of the interfering source at electrical angle $\theta_2=\pi/n$ (half the Rayleigh limit of the array). For our simulations, we use Gaussian compression matrices $\mathbf{\Phi}_{m\times n}$ whose elements are i.i.d. $\mathcal{CN}(0,1/m)$.  The Fisher information matrix and the CRB on the estimation of $\theta_1$ are calculated for different realizations of $\Phib$. Fig. \ref{fig:Hist} shows the CRB on the estimation of $\theta_1$ before compression divided by its corresponding value after compression, i.e. $\frac{(\Jb^{-1})_{11}}{(\Jbh^{-1})_{11}}$ for $m=64$, $n=128$. A histogram of actual values of $\frac{(\Jb^{-1})_{11}}{(\Jbh^{-1})_{11}}$ for $10^5$ independent realizations of random $\mathbf{\Phi}$ is shown in blue. The red curve represents the pdf of a $B^I(m-p+1,n-m)$ distributed random variable for $p=2$. This figure simply provides an illustration of the result \eqref{eq:normJ}. 

\begin{figure}[ht]\centering
  \includegraphics[width=288pt]{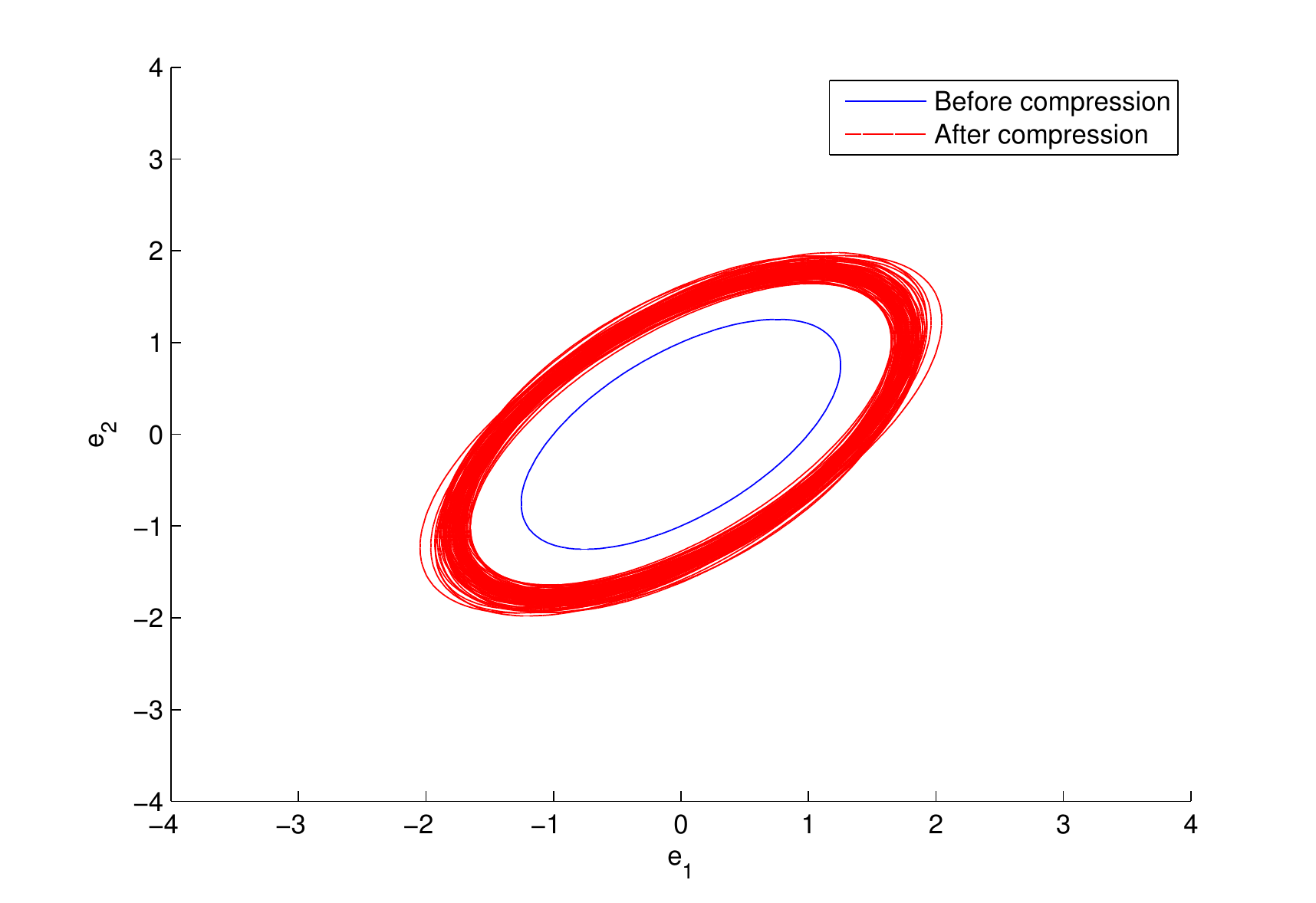}
  \caption{Concentration ellipses for the Fisher information matrices before and after compression.}\label{fig:CE}
\end{figure}

Recall that the inverse Fisher Information matrix $\Jb^{-1}$ lower bounds the error covariance matrix $\bm{\Sigma}=E[\eb\eb^H]$ for unbiased errors $\eb=\hat{\bm{\theta}}-\bm{\theta}$. So the concentration ellipse $\eb^{H}\bm{\Sigma}^{-1}\eb\le \eb^H\Jb\eb$ for all $\eb\in\mathbb{C}^{p}$. The ellipses $\eb^H\Jb\eb=r^2$ and $\eb^H\Jbh\eb$, with $r^2=\Jb_{11}$, are illustrated in Fig. \ref{fig:CE}, demonstrating the effect that compression inflates the concentration ellipse. The blue curve is the locus of all points $\mathbf{e}\in\mathbb{C}^{p}$, for which $\eb^H\Jb\eb=r^2$. The red curves are the loci of all points $\eb\in\mathbb{C}^{p}$, for which $\eb^H\Jbh\eb=r^2$ for $100$ realizations of the Fisher information matrix after compression. As can be seen, the concentration ellipse for the Fisher information matrix before compression has the smallest volume in comparison with all the realizations of the concentration ellipses after compression. Also, for each realization of the Gaussian compression, the orientation of the concentration ellipse is nearly aligned with that of the uncompressed ellipse.   

\begin{figure}[ht]\centering
  \includegraphics[width=288pt]{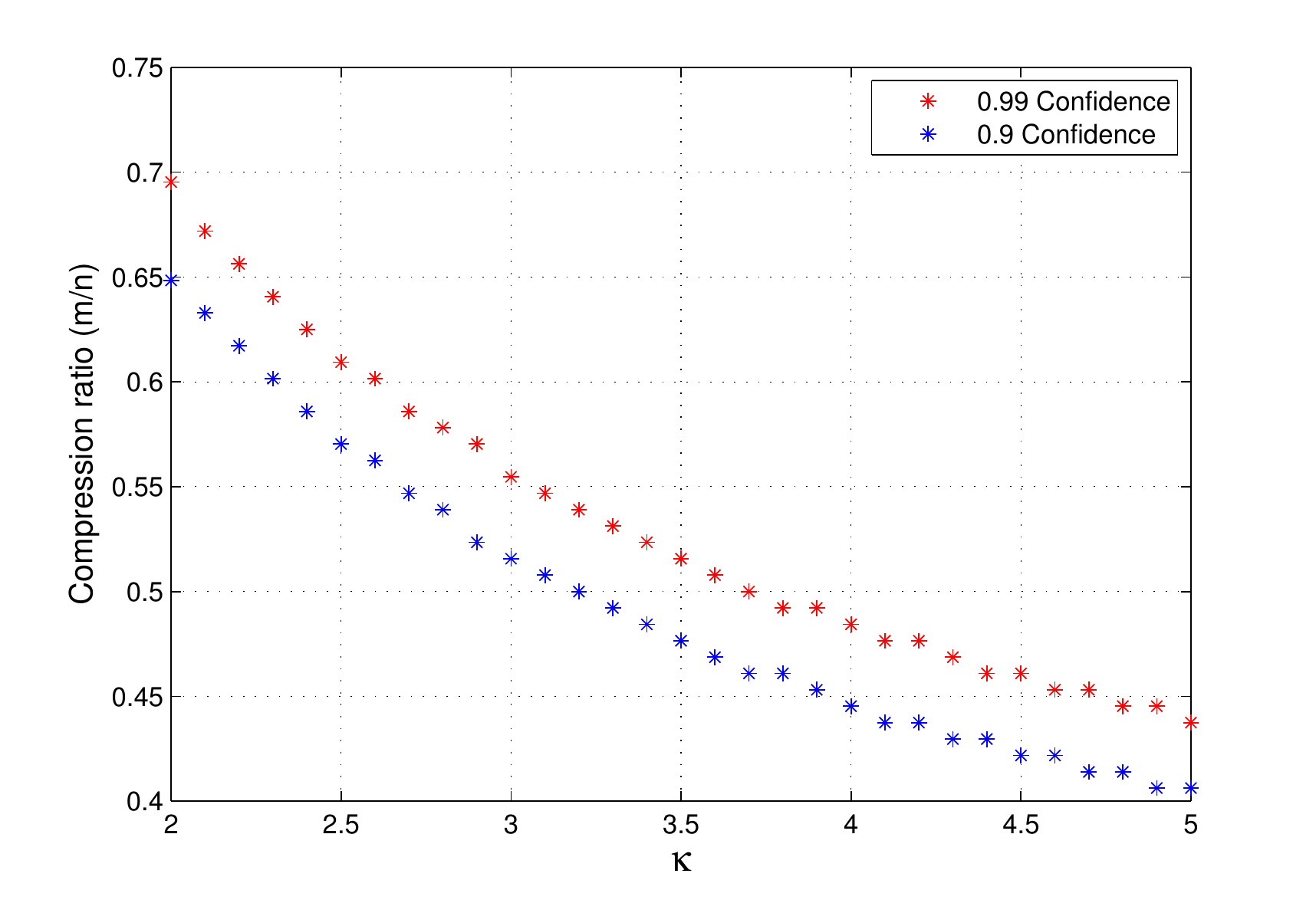}
  \caption{Compression ratios needed so that $(\Jbh^{-1})_{11}<\kappa(\Jb^{-1})_{11}$ for different confidence levels.}\label{fig:Prob}
\end{figure}

Figure \ref{fig:Prob} shows the compression ratio $m/n$ needed so that the CRB after compression $(\Jbh^{-1})_{11}$ does not exceed $\kappa$ times the CRB before compression $(\Jb^{-1})_{11}$, at two levels of confidence and for $n=128$. These curves are plotted using the tail probabilities of a univariate beta random variables. They can be used as guidelines for deriving a satisfactory compression ratio based on a tolerable level of loss in the CRB. Alternatively, we can plot the confidence level curves versus $m$ for fixed values of $\kappa$ . In that case, the plots may be useful to find a number of measurements that would guarantee that after compression CRB does not go above a desired bound (corresponding to a particular $\kappa$) with a certain level of confidence.

\section{Conclusion}
In this paper, we have studied the effect of random compression of noisy measurements on the CRB for estimating parameters in a nonlinear model. We have considered the class of random compression matrices whose distributions are right-orthogonally invariant. A random compression matrix with i.i.d. standard normal elements is one such compression matrix. The analytical distributions obtained in this paper can be used to quantify the amount of loss due to compression. Also, they can be used as guidelines for choosing a suitable compression ratio based on a tolerable loss in the CRB. Importantly, the distribution of the ratios of CRBs before and after compression depends only on the number of parameters and the number of measurements. The distribution is invariant to the underlying signal-plus-noise model, in the sense that it is invariant to the underlying (before compression) Fisher information matrix.

\bibliographystyle{IEEEtran}
\bibliography{ref_CRB}

\begin{thebibliography}{10}
\providecommand{\url}[1]{#1}
\csname url@samestyle\endcsname
\providecommand{\newblock}{\relax}
\providecommand{\bibinfo}[2]{#2}
\providecommand{\BIBentrySTDinterwordspacing}{\spaceskip=0pt\relax}
\providecommand{\BIBentryALTinterwordstretchfactor}{4}
\providecommand{\BIBentryALTinterwordspacing}{\spaceskip=\fontdimen2\font plus
\BIBentryALTinterwordstretchfactor\fontdimen3\font minus
  \fontdimen4\font\relax}
\providecommand{\BIBforeignlanguage}[2]{{%
\expandafter\ifx\csname l@#1\endcsname\relax
\typeout{** WARNING: IEEEtran.bst: No hyphenation pattern has been}%
\typeout{** loaded for the language `#1'. Using the pattern for}%
\typeout{** the default language instead.}%
\else
\language=\csname l@#1\endcsname
\fi
#2}}
\providecommand{\BIBdecl}{\relax}
\BIBdecl

\bibitem{CS1}
E.~J. Cand\`{e}s, ``Compressive sampling,'' in \emph{Proc. Int. Congress
  Math.}, vol.~3, 2006, pp. 1433--1452.

\bibitem{CS2}
D.~L. Donoho, ``Compressed sensing,'' \emph{IEEE Transactions on Information
  Theory}, vol.~52, no.~4, pp. 1289--1306, 2006.

\bibitem{CS3}
R.~Baraniuk, ``Compressive sensing,'' \emph{IEEE Signal Processing Magazine},
  vol.~24, no.~4, pp. 118--121, Jul. 2007.

\bibitem{chi-sp11}
Y.~Chi, L.~Scharf, A.~Pezeshki, and A.~Calderbank, ``Sensitivity to basis
  mismatch in compressed sensing,'' \emph{IEEE Transactions on Signal
  Processing}, vol.~59, no.~5, pp. 2182--2195, May 2011.

\bibitem{chi-icassp10}
Y.~Chi, A.~Pezeshki, L.~L. Scharf, and R.~Calderbank, ``Sensitivity to basis
  mismatch in compressed sensing,'' in \emph{Proc. 2010 IEEE Int. Conf. on
  Acoust., Speech and Signal Process. (ICASSP), Dallas, TX}, Mar. 2010, pp.
  3930--3933.

\bibitem{scharf-asilomar11}
L.~L. Scharf, E.~K.~P. Chong, A.~Pezeshki, and J.~R. Luo, ``Sensitivity
  considerations in compressed sensing,'' in \emph{Conf. Rec. 45th Annual
  Asilomar Conf. Signals, Systs., Computs., Pacific Grove, CA,}, Nov. 2011, pp.
  744--748.

\bibitem{scharf-DASP11}
------, ``Compressive sensing and sparse inversion in signal processing:
  Cautionary notes,'' in \emph{Proc. 7th Workshop on Defence Applications of
  Signal Processing (DASP), Coolum, Queensland, Australia, Jul. 10-14, 2011}.

\bibitem{Pakrooh13a}
P.~Pakrooh, L.~L.~Scharf, A.~Pezeshki, and Y.~Chi, ``Analysis of {F}isher
  information and the {C}ram\'{e}r-{R}ao bound for nonlinear parameter
  estimation after compressed sensing,'' in \emph{Proc. IEEE International
  Conference on Acoustics, Speech and Signal Processing (ICASSP)}, Vancouver,
  BC, Canada, May 2013, pp. 6630--6634.

\bibitem{Babadi09}
B.~Babadi, N.~Kalouptsidis, and V.~Tarokh, ``Asymptotic achievability of the
  {C}ram\'{e}r-{R}ao bound for noisy compressive sampling,'' \emph{IEEE
  Transactions on Signal Processing}, vol.~57, no.~3, pp. 1233--1236, Mar.
  2009.

\bibitem{Babaizadeh12}
R.~Niazadeh, M.~Babaie-Zadeh, and C.~Jutten, ``On the achievability of
  {C}ram\'{e}r-{R}ao bound in noisy compressed sensing,'' \emph{IEEE
  Transactions on Signal Processing}, vol.~60, no.~1, pp. 518--526, Jan. 2012.

\bibitem{Madhow12}
D.~Ramasamy, S.~Venkateswaran, and U.~Madhow, ``Compressive estimation in
  {A}{W}{G}{N}: {G}eneral observations and a case study,'' in \emph{Conf. Rec.
  46th Annual Asilomar Conf. Signals, Systs., Computs., Pacific Grove, CA,},
  Nov. 2012, pp. 953--957.

\bibitem{Nielsen2012}
J.~K. Nielsen, M.~G. Christensen, and S.~H. Jensen, ``On compressed sensing and
  the estimation of continuous parameters from noisy observations,'' in
  \emph{Proc. IEEE International Conference on Acoustics, Speech and Signal
  Processing (ICASSP)}, Kyoto, Japan, Mar. 2012, pp. 3609--3612.

\bibitem{4-Eldar}
S.~M. Kay, \emph{Fundamentals of Statistical Signal Processing, Volume I:
  Estimation Theory}.\hskip 1em plus 0.5em minus 0.4em\relax Prentice-Hall,
  Apr. 1993.

\bibitem{Scharf_geo}
L.~L. Scharf and L.~T. McWhorter, ``Geometry of the {C}ram\'{e}r-{R}ao bound,''
  \emph{Signal Process.}, vol.~31, no.~3, pp. 301--311, Apr. 1993.

\bibitem{Srivastava}
M.~Srivastava, ``On the complex {W}ishart distribution,'' \emph{The Annals of
  Mathematical Statistics}, pp. 313--315, 1965.

\bibitem{RMB}
I.~S. Reed, J.~D. Mallett, and L.~E. Brennan, ``Rapid convergence rate in
  adaptive arrays,'' \emph{IEEE Transactions on Aerospace and Electronic
  Systems}, no.~6, pp. 853--863, 1974.

\bibitem{James1964}
A.~T. James, ``Distributions of matrix variates and latent roots derived from
  normal samples,'' \emph{The Annals of Mathematical Statistics}, pp. 475--501,
  1964.

\bibitem{Chikuse2003}
Y.~Chikuse, \emph{Statistics on special manifolds}.\hskip 1em plus 0.5em minus
  0.4em\relax New York: Springer, 2003.

\bibitem{Gupta2005}
A.~K. Gupta, D.~K. Nagar, and E.~Bedoya, ``Properties of the complex matrix
  variate {D}irichlet distribution,'' \emph{Scientiae Mathematicae Japonicae},
  vol.~66, no.~1, p.~53, 2007.

\end{thebibliography}

\end{document}